\input amstex
\documentstyle{amsppt}
\magnification=\magstep1 

\topmatter 
\title Some observations on compact indestructible spaces
\endtitle 
\author Angelo Bella \endauthor 
\affil University of Catania, Italy \endaffil 

\address
{Department of Mathematics -- University of Catania \newline
Citta' universitaria  viale A. Doria 6 \newline
95125 Catania, Italy}
\endaddress

\email
bella\@dmi.unict.it
\endemail

\abstract {Inspired by a recent work of Dias and Tall, 
we show that a compact indestructible space is
sequentially compact.  We also prove that a Lindel\" of $T_2$ 
indestructible space has the finite derived set property  and a
compact  
$T_2$ indestructible space is
pseudoradial. } 
\endabstract 
\endkeywords 
\endsubjclass 
\endtopmatter 
\document
 A compact space is indestructible if it remains compact in any
countably closed forcing extension. 
This is a particular case  of the notion of 
Lindel\"of indestructibility,    whose study was initiated by
Tall in
\cite{9}. A space is compact indestructible if and only if it is
compact and Lindel\"of indestrutible. 
  A nice connection of Lindel\"of indestructibility    with 
certain infinite topological game was later 
 discovered by Scheepers and Tall \cite{8}.

   ${\text G}^{\omega_1}
_1(\Cal O,
\Cal O)$  denotes  the
game of length $\omega_1 $ played on a topological space $X$ by
two players I and II     in the
following way:   at  the $\alpha $-th inning  
player I choose an open cover $\Cal U_\alpha $ of $X$  and
player II  responds
by taking an element $U_\alpha  \in \Cal U_\alpha $. Player II
wins if and
only if   $\{U_\alpha 
:\alpha <\omega_1 \} $ covers $X$. 
\proclaim {Proposition 1}(\cite{8}, Theorem 1) A Lindel\"of
space $X$ is indestructibily Lindel\" of if and only if player I
does not have a winning strategy in ${\text G}_1^{\omega_1}(\Cal
O, \Cal
O)$. \endproclaim   

   Recently,  Dias and Tall \cite{4} started to  
investigate the topological structure of compact indestructible
spaces. In particular, they proved that a compact $T_2$
indestructible space contains a non-trivial convergent sequence
(\cite{4}, Corollary 3.4). 

The aim of this short note is to strengthen the above result, by
showing that indestructibility  actually gives even more than
sequential
compactness (Theorem 3).  However, indestructibility forces a
compact space to be sequentially compact  in the absolute general
case, that is by  assuming no 
separation axiom (Theorem 1).   
The same proof, with minor changes, will show
that a Lindel\" of $T_2$ indestructible space has the  finite
derived set property
(Theorem 2).

As usual, $A\subseteq ^* B$ means  $|A\setminus B|<\aleph_0$ (mod
finite inclusion).
\proclaim {Theorem 1} Every compact indestructible space is
sequentially compact. \endproclaim 
\demo {Proof} Let $X$ be a compact indestructible space and
assume that $X$ is not sequentially compact. Our task is to show
that in this case player I  would have a winning strategy in the
game
${\text G}_1^{\omega_1}(\Cal O, \Cal O)$. Fix a countable
infinite set
$A\subseteq X$ with no  infinite convergent subsequence. For
each $x\in X$ there is an open set $U_x$ such that $x\in U_x$
and $|A\setminus U_x|=\aleph_0$.  The first move of player I is
the open cover $\Cal U_0=\{U_x:x\in X\}$.  If player II responds
by choosing $U_{\{x_0\}}\in \Cal U_0$, then let
$A_{\{x_0\}}=A\setminus U_{\{x_0\}}$. For each $x\in X$ there is
an open set $U_{(x_0,x)}$ such that $x\in U_{(x_0,x)}$ and
$|A_{\{x_0\}}\setminus U_{(x_0,x)}|=\aleph_0$.   The second move
of player I  is  the open cover $\Cal U_1=\{U_{(x_0,x)} :x\in
X\}$.  If 
 player II responds by choosing $U_{(x_0,x_1)}\in \Cal U_1$, 
then  let $A_{(x_0,x_1)}=A_{\{x_0\}}\setminus U_{(x_0,x_1)}$.
Again, for each $x\in X$ player I chooses an open set
$U_{(x_0,x_1,x)}$ such that $x\in U_{(x_0,x_1,x)}$ and
$|A_{(x_0,x_1)}\setminus U_{(x_0,x_1,x)}|=\aleph_0$. At
the $\omega$-th inning of the game,  the moves of the two players
have  defined a function
$f:\omega \to X$ and a decreasing chain of sets
$\{A_{f\restriction n} :n<\omega\}$.  Player I  chooses an
infinite 
set $B\subseteq A$  satisfying $B\subseteq ^* A_{f\restriction
n}$ for each $n<\omega$ and for each $x\in X$ an open set
$U_{f\frown x}$ such that $x\in U_{f\frown x}$ and $|B\setminus
U_{f\frown x}|=\aleph_0$. Then, at the $\omega$-th inning player
I plays the open cover $\Cal U_\omega=\{U_{f\frown x} : x\in
X\}$. If player II responds by choosing $U_{f\frown x}$, then let
$x_\omega=x$ and $A_f=B\setminus U_{f\frown x_\omega}$. In
general, at the $\alpha $-th inning the moves of the two
players have already defined a function $f:\alpha \to X$ and a
mod
finite decreasing family $\{A_{f\restriction \beta}: \beta<\alpha
\}$ of infinite subsets of  $A$.  Then, player I fixes an
infinite
set $B\subseteq A$  such that $B\subseteq ^* A_{f\restriction
\beta}$ for each $\beta<\alpha $ and plays the open cover $\Cal
U_\alpha
=\{U_{f\frown x}: x\in X\}$, where  $x\in U_{f\frown x}$ and
$|B\setminus U_{f\frown x}|=\aleph_0$. If the responds of player
II is $U_{f\frown x }$, then let $x_\alpha =x$, $A_f=B\setminus
U_{f\frown x_\alpha }$
and so on.  

 At the end of the game, we  
have  a
function $g:\omega_1\to X$ and  a mod finite decreasing chain
$\{A_{g\restriction \alpha } :\alpha <\omega_1\}$ of infinite
subsets of $A$.   The set resulting  from the moves of player II
is
the collection $\Cal V=\{U_{g\restriction {\alpha +1}}: \alpha
<\omega_1\}$. For any finite
set of ordinals $\alpha _0, \ldots, \alpha _m<\omega_1$, taking
some $\beta<\omega_1$ such that $\alpha _i<\beta$ for  $i\le m$, 
we see that  the infinite set $A_{g\restriction \beta}$  has a
finite intersection with each $U_{g\restriction {\alpha _i+1}}$
and therefore  the
subcollection $\{U_{g\restriction {\alpha _i+1}}: i\le m\}$
cannot cover $X$.  Since 
$\Cal V$  does not have finite subcovers, the compactness of $X$
implies that the whole $\Cal V$ cannot cover $X$. Thus,   player
I wins the game, in contrast with  
Proposition 1.  \qed\enddemo 

Recall that a topological space $X$ has the {\it finite derived
set}  (briefly FDS) property   provided that every infinite set
of
$X$ contains an infinite subset with  at most finitely many
accumulation points (see for instance \cite{2}).   Since in a 
$T_2$ space a convergent sequence has only one accumulation
point, we see that if a     $T_2$ space has a countable infinite
set $A$ violating  the  finite derived set property, then for
each infinite set $B\subseteq A$ and  each point $x\in X$ there
must be an open set $U_x$ such that $x\in U_x$ and $|B\setminus
U_x|=\aleph_0$. Notice, however, that for this much less than
$T_2$ is needed. For instance, it suffices  for the space  to be
SC, namely  that 
every convergent sequence
together with the limit point is a closed subset 
(see \cite{2}).  

 With this observation in mind, we can 
modify the above proof to get the following :
\proclaim {Theorem 2} A Lindel\"of $T_2$ indestructible space has
the finite derived set property. \endproclaim
\demo {Proof}
 Let $X$ be a  Lindel\" of $T_2$ indestructible space and
assume that $X$ does not have the FDS property. As in the proof
of Theorem 1,   our task is to show
that in this case player I would have a winning strategy in the
game
${\text G}_1^{\omega_1}(\Cal O, \Cal O)$. Fix a countable
infinite set
$A\subseteq X$ witnessing the failure of the FDS property.
Taking into account the paragraph before the theorem,     for
each infinite set $B\subseteq A$ and each $x\in X$ there is an
open set $U_x$ such that $x\in U_x$
and $|B\setminus U_x|=\aleph_0$. Now, the strategy of player I is
exactly the same as that  in the proof of Theorem 1.   At the end
of
the game,    The set resulting  from the moves of player II is
again the collection $\Cal V=\{U_{g\restriction {\alpha +1}}:
\alpha
<\omega_1\}$. We claim that $\Cal V$ cannot cover $X$. Otherwise,
by the Lindel\"ofness of $X$,  there should exists a countable
set of ordinals $S\subseteq \omega_1$ such that  the
subcollection  $\{U_{g\restriction {\alpha +1}}:\alpha \in S\}$
 would cover $X$.  Taking
some $\beta<\omega_1$ such that $\alpha <\beta$ for each $\alpha
\in S$,  
we see that  the infinite set $A_{g\restriction \beta}$  has a
finite intersection with $U_{g\restriction {\alpha +1}}$ for each
$\alpha \in S$. But, this implies  that the infinite set
$A_{g\restriction \beta}$ does not have accumulation points in
$X$, in contrast with the supposed failure of the FDS property  
in $A$. Thus,   $\Cal V$ cannot cover $X$ and again  player
I wins the game.  \qed\enddemo

The above theorem   provides  new informations on the topological
structure of a Lindel\" of indestrutible space.

 We will finish by showing that  for $T_2$ spaces Theorem 1 can
be improved.
\proclaim  {Proposition 2} (\cite{4}, Corollary 3.3)  A compact
$T_2$ space  which is not
first countable at any point is destructible. \endproclaim 
Recall that a topological space $X$ is pseudoradial   provided
that for any non-closed set $A\subseteq X$ there exists a 
well-ordered net  $S\subseteq A$ which converges to a point
outside
$A$. For more on these spaces see \cite{3}.

Clearly every compact pseudoradial space is sequentially compact,
but 
 the converse may consistently fail \cite{5}.
\proclaim {Theorem 3} Any compact $T_2$ indestructible space is
pseudoradial. \endproclaim 
\demo {Proof} Let $X$ be a compact $T_2$ indestructible space and
let $A$ be a non-closed subset. Let $\lambda $ be the smallest
cardinal such that there exists a non-empty $G_\lambda $-set
$H\subseteq \overline A\setminus A$. As $X$ is indestructible, so
is the subspace $H$. Hence, by  Proposition 2, 
$H$ is first countable at some point $p $.
Clearly, $\{p\}$ is a $G_\lambda $-set in $X$ and so there are
open sets $\{U_\alpha :\alpha <\lambda \}$ satisfying
$\{p\}=\bigcap \{U_\alpha :\alpha <\lambda \}=\bigcap \{\overline
{U_\alpha }:\alpha <\lambda \}$. The minimality of $\lambda $
ensures that for each $\alpha <\lambda $ we may pick a point
$x_\alpha \in A\cap \bigcap\{U_\beta:\beta<\alpha \}$. The
compactness of $X$ implies that the  well-ordered net $\{x_\alpha
:\alpha
<\lambda \}$ converges to $p$ and we are done. \qed\enddemo 

Notice that the indestructibility of a compact space is stronger
than pseudoradiality: the Example in section 3 of  \cite{4} 
is a compact
$T_2$ pseudoradial space which is destructible. 

 The fact that   pseudoradiality is a weakening
of sequentiality and the well-know       
fact 
that   compact spaces  of 
countable
tightness  are sequential  under PFA \cite{1}, might suggest
that a
compact $T_2$ indestructible space of countable tightness is
always 
sequential. But this is not the case:   the one-point
compactification
of the Ostaszewski's space \cite{7} is a non-sequential compact
$T_2$ space
of countable tightness  which is  indestructible having
cardinality $\aleph_1$ (see \cite{4}).

Theorem 3 is no longer true for Lindel\" of spaces.  Koszmider
and Tall constructed
\cite{6} a
model of ZFC+CH where there exists a regular Lindel\"of $P$-space
$X$ of cardinality $\aleph_2$ without Lindel\"of subspaces of
size
$\aleph_1$. Such a space does not  have 
convergent well-ordered nets of length $\aleph_1$.  Therefore,
$X$ is not pseudoradial because it obviously  contains
non-closed subsets of cardinality
$\aleph_1$.  On the other hand, it is easy to check that a
Lindel\"of  
P-space is indestructible (see e. g. \cite{8}).

Recall that a space $X$ satisfies the selection principle
${\text S}_1^{\omega_1}(\Cal O, \Cal O)$ provided that for any 
family 
$\{\Cal U_\alpha : \alpha <\omega_1\}$ of open covers of $X$ one
may pick an
element $U_\alpha \in \Cal U_\alpha $ in such
a way that the
collection $\{U_\alpha :\alpha <\omega_1\}$ covers $X$.

It is clear that  any compact indestructible space satisfies
${\text S}_1^{\omega_1}(\Cal O, \Cal O)$ and the example
described in
section 3 of \cite{4} shows that the previous
implication is consistently not reversible.  Such example is a
compact LOTS
and
so it is sequentially compact. An obvious question then arises:
\proclaim {Question 1} Let $X$  be a compact (or compact $T_2$)
space satisfying ${\text S}_1^{\omega_1}(\Cal O, \Cal O)$. Is $X$
sequentially compact? 
\endproclaim 

A space answering the above question in the negative   would
provide
a compact space satisfying 
${\text S}_1^{\omega_1}(\Cal O, \Cal O)$ which is ``more
destructible''
than the mentioned example in \cite{4}. 

We can also formulate a weaker version of the problem.
\proclaim {Question 2} Is it true that  any compact $T_2$ space
satisfying ${\text S}_1^{\omega_1}(\Cal O, \Cal O)$ contains a
non-
trivial convergent sequence? \endproclaim 
An interesting feature of the above question is that  any
counterexample to it turns out to be  an Efimov's space, that is 
a
compact
$T_2$ space containing no copy of $\beta\omega$ and no 
non-trivial convergent sequence.

\enddocument